\documentclass[12pt,psamsfonts]{amsart}

\newtheorem{anyprop}{Anyprop}[section]

\newtheorem{theorem}[anyprop]{Theorem}

\newtheorem{proposition}[anyprop]{Proposition}
\newtheorem{corollary}[anyprop]{Corollary}

\newtheorem{theoremintro}{Theorem}

\theoremstyle{definition}

\newtheorem{situation}[anyprop]{Situation}

\newtheorem{example}[anyprop]{Example}

\newtheorem{remark}[anyprop]{Remark}

\renewcommand{\AA}{\mathbb{A}}

\newcommand{\ZZ}{\mathbb{Z}}
\newcommand{\QQ}{\mathbb{Q}}

\newcommand{\PP}{\mathbb{P}}

\newcommand  {\shG}     {\mathcal{G}}

\newcommand  {\shI}     {\mathcal{I}}

\newcommand  {\shL}     {\mathcal{L}}

\newcommand  {\shQ}     {\mathcal{Q}}

\newcommand  {\shS}     {\mathcal{S}}

\newcommand  {\fom}     {\mathfrak{m}}

\newcommand  {\engrad}     {\operatorname{end}}

\newcommand  {\Hom}     {\operatorname{Hom}}

\newcommand  {\im}      {\operatorname{im}}

\renewcommand  {\ker }  {\operatorname{kern}}

\newcommand  {\lra}     {\longrightarrow}

\renewcommand{\O}       {\mathcal{O}}

\newcommand  {\Proj}    {\operatorname{Proj}}

\newcommand  {\ra}      {\rightarrow}

\newcommand  {\rk}    {\operatorname{rk}}

\newcommand  {\reg}     {\operatorname{reg}}

\newcommand  {\Spec}    {\operatorname{Spec}}

\newcommand  {\Syz}     {\operatorname{Syz}}

\newcommand{\komdots}{ , \ldots , }
\newcommand{\plusdots}{ + \ldots + }
\newcommand{\minusdots}{ - \ldots - }

\newcommand{\congdots}{ \cong \ldots \cong }

\newcommand{\lradots}{ \lra \ldots \lra }

\newcommand{\subsetdots}{ \subset \ldots \subset }

\theoremstyle{remark}

\numberwithin{equation}{section}

\usepackage{amscd}
\usepackage{amssymb}

\def\mydate{\number\day\space\ifcase\month \or January\or February\or March\or April\or May\or
June\or July\or August\or September\or October\or November\or
December\fi \space\number\year}

\newcommand{\dimty}{t}

\begin{document}

\title[Bounds for Frobenius powers and for tight closure]
{A linear bound for Frobenius powers and an inclusion bound for tight closure}

\author[Holger Brenner]{Holger Brenner}
\address{Department of Pure Mathematics, University of Sheffield,
  Hicks Building, Hounsfield Road, Sheffield S3 7RH, United Kingdom}
\email{H.Brenner@sheffield.ac.uk}





\begin{abstract}
Let $I$ denote an $R_+$-primary homogeneous ideal in a normal
standard-graded Cohen-Macaulay domain over a field of positive
characteristic $p$. We give a linear degree bound for the Frobenius
powers $I^{[q]}$ of $I$, $q=p^{e}$, in terms of the minimal slope of
the top-dimensional syzygy bundle on the projective variety $\Proj
R$. This provides an inclusion bound for tight closure. In the same
manner we give a linear bound for the Castelnuovo-Mumford regularity
of the Frobenius powers $I^{[q]}$.
\end{abstract}

\maketitle


\section*{Introduction}

Let $R$ denote a noetherian ring, let $\fom$ denote a maximal ideal in $R$ and let $I$
denote an $\fom$-primary ideal. This means by definition that $\fom$ is the radical of $I$.
Then there exists a (minimal) number $k$ such that
$\fom^k \subseteq I \subseteq \fom$ holds.
If $R$ contains a field of positive characteristic $p$,
then the Frobenius powers of the ideal $I$, that is
$$I^{[q]} = \{ f^q: f \in I\} \,  , \, \, \, q = p^{e} ,$$
are also $\fom$-primary and hence there exists a minimal number
$k(q)$ such that $\fom^{k(q)} \subseteq I^{[q]}$ holds. In this
paper we deal with the question how $k(q)$ behaves as a function of
$q$, in particular we look for linear bounds for $k(q)$ from above.
If $\fom^k \subseteq I$ and if $l$ denotes the number of generators
for $\fom^k$, then we get the trivial linear inclusion
$(\fom^k)^{lq} \subseteq (\fom^k)^{[q]} \subseteq I^{[q]}$.

The main motivation for this question comes from the theory of tight
closure. Recall that the tight closure of an ideal $I$ in a domain
$R$ containing a field of positive characteristic $p$ is the ideal
$$I^* = \{f \in R: \exists  0 \neq c \in R \mbox{ such that }
cf^q \in I^{[q]} \mbox{ for all } q=p^{e} \} \, .$$ A linear
inclusion relation $ \fom^{\lambda q + \gamma} \subseteq I^{[q]}$
for all $q=p^{e}$ implies the inclusion $\fom^\lambda \subseteq
I^*$, since then we can take any element $0 \neq c \in \fom^\gamma$
to show for $f \in \fom^\lambda$ that $cf^q \in \fom^{\lambda q
+\gamma} \subseteq I^{[q]}$, hence $f \in I^*$. The trivial bound
mentioned above yields $\fom^{kl} \subseteq I^*$, but in fact we
have already $ \fom^{kl} \subseteq \fom^k\subseteq I$, so this does
not yield anything interesting.

We restrict in this paper to the case of a normal standard-graded
domain $R$ over an algebraically closed field $K=R_0$ of positive
characteristic $p$ and a homogeneous $R_+$-primary ideal $I$. The
question is then to find the minimal degree $k(q)$ such that
$R_{\geq k(q)} \subseteq I^{[q]}$ or at least a good linear bound
$k(q) \leq \lambda q+\gamma$. In this setting we work mainly over
the normal projective variety $Y=  \Proj R$, endowed with the very
ample invertible sheaf $\O_Y(1)$. If $I=(f_1 \komdots f_n)$ is given
by homogeneous ideal generators $f_i$ of degree $d_i = \deg (f_i)$,
then we get on $Y$ the following short exact sequences of locally
free sheaves,
$$0 \lra \Syz(f_1^q \komdots f_n^q)(m) \lra \bigoplus_{i=1}^n \O_Y(m -qd_i)
\stackrel{f_1^q \komdots f_n^q}{\lra} \O_Y(m) \lra 0 \, .$$ Another
homogeneous element $h \in R$ of degree $m$ yields a cohomology
class $\delta(h) \in H^1(Y,\Syz(f_1^q \komdots f_n^q)(m))$, and
therefore the question whether $h \in (f_1^q \komdots
f_n^q)=I^{[q]}$ is equivalent to the question whether $\delta(h)=0$.
Since $\Syz(f_1^q \komdots f_n^q)(0) =F^{*e} (\Syz(f_1 \komdots
f_n)(0))$ is the pull-back under the $e$-th absolute Frobenius
morphism $F^{e}:Y \ra Y$, our question is an instance of the
following more general question: given a locally free sheaf $\shS$
on a normal projective variety $(Y,\O_Y(1))$, find an
(affine-linear) bound $\ell (q)$ such that for $m \geq \ell (q)$ we
have $H^1(Y, \shS^q (m))=0$, where we set $\shS^{q}=F^{e*}(\shS)$.
Using a resolution $\shG_\bullet \ra \shS \ra 0$ where $\shG_j=
\bigoplus_{(k,j)} \O_Y( -\alpha_{k,j})$, we can shift the problem
(at least if $Y= \Proj R$ with $R$ Cohen-Macaulay, so that
$H^{i}(Y,\O_Y(m))=0$ for $0< i < \dim (Y)$) to the problem of
finding a bound such that $H^t(Y, \shS_t^q(m)) =0$, where $\shS_t =
\ker (\shG_t \ra \shG_{t-1})$ and $t = \dim (Y)$. By Serre duality
this translates to $\Hom( \shS_t^q(m), \omega_Y)=0$. Now the
existence of such mappings is controlled by the minimal slope of
$\shS_t^q(m)$. Let $\bar{\mu}_{\min} (\shS_t)= \lim \inf_{q=p^{e}}
\mu_{\min}(\shS_t^q)/q$ and set $\nu =- \bar{\mu}_{\min} (\shS_t)/
\deg (Y)$. With these notations applied to $\shS=\Syz(f_1 \komdots
f_n)(0)$ our main results are the following theorems (Theorems
\ref{theoreminclusion} and \ref{tightinclusion}).

\begin{theoremintro}
\label{theoreminclusionintro}
Let $R$ denote a standard-graded normal Cohen-Macaulay domain over
an algebraically closed field $K$ of characteristic $p >0$.
Suppose that the dualizing sheaf $\omega_Y$ of $Y=  \Proj R$ is invertible.
Let $I$ denote a homogeneous $R_+$-primary ideal.
Then $R_{ > q \nu + \frac{\deg (\omega_Y)}{\deg(Y)} } \subseteq I^{[q]}$.
\end{theoremintro}

From this linear bound for the Frobenius powers
we get the following inclusion bound for tight closure.

\begin{theoremintro}
\label{tightinclusionintro} Under the assumptions of Theorem
\ref{theoreminclusionintro} we have the inclusion $R_{\geq  \nu}
\subseteq I^*$, where $I^*$ denotes the tight closure of $I$.
\end{theoremintro}

This theorem generalizes \cite[Theorem 6.4]{brennerslope} from
dimension two to higher dimensions. We also obtain an inclusion
bound for the Frobenius closure (Corollary \ref{frobenius}) and a
linear bound for the Castelnuovo-Mumford regularity of the Frobenius
powers $I^{[q]}$ (Theorem \ref{regularitybound}), which improves a
recent result of M. Chardin \cite{chardinregularitypowers}.

I thank M. Blickle for useful remarks.

\section{Some projective preliminaries}

Let $K$ denote an algebraically closed field and let $Y$ denote a
normal projective variety over $K$ of dimension $t$ together with a
fixed ample Cartier divisor $H$ with corresponding ample invertible
sheaf $\O_Y(1)$. The degree of a coherent torsion-free sheaf $\shS$
(with respect to $H$) is defined by the intersection number
$\deg(\shS)= \deg (c_1(\shS) )= c_1(\shS) . H^{t-1}$, see
\cite[Preliminaries]{maruyamagrauertmuelich} for background of this
notion. The degree is additive on short exact sequences \cite[Lemma
1.5(2)]{maruyamagrauertmuelich}.

The slope of $\shS$ (with respect to $H$), written $\mu(\shS)$, is
defined by dividing the degree through the rank. The slope fulfills
the property that $\mu (\shS_1 \otimes \shS_2)= \mu(\shS_1) +
\mu(\shS_2)$ \cite[Lemma 1.5(4)]{maruyamagrauertmuelich}. The
minimal slope of $\shS$, $\mu_{\min} (\shS)$, is given by
$$\mu_{\min} (\shS) = \inf \{ \mu(\shQ):\, \shS \ra \shQ \ra 0
\mbox{ is a torsion-free quotient sheaf} \} \, .$$ If $\shS_1
\subsetdots \shS_k= \shS$ is the Harder-Narasimhan filtration of
$\shS$ \cite[Proposition 1.13]{maruyamagrauertmuelich}, then
$\mu_{\min} (\shS) = \mu(\shS/\shS_{k-1})$. If $\shL$ is an
invertible sheaf and $\mu_{\min}(\shS) > \deg(\shL)$, then there
does not exist any non-trivial sheaf homomorphism $\shS \ra \shL$.
The sheaf $\shS$ is called semistable if $ \mu(\shS)
=\mu_{\min}(\shS)$.

Suppose now that the characteristic of $K$ is positive and let
$F^{e}\!:Y \ra Y$ denote the $e$-th absolute Frobenius morphism. We
denote the pull-back of $\shS$ under this morphism by $\shS^q
=F^{e*}(\shS)$, $q=p^{e}$. The slope behaves like $\mu(\shS^q) = q
\mu (\shS)$ (this follows from \cite[Lemma
1.6]{maruyamagrauertmuelich}, for which it is enough to assume that
the finite mapping is flat in codimension one; note that we compute
the slope always with respect to $\O_Y(1)$, not with respect to
$F^{*e}(\O_Y(1))=\O_Y(q)$). It may however happen that $\mu_{\min}
(\shS^q) < q \mu_{\min} (\shS)$. Therefore it is useful to consider
the number (compare \cite{langersemistable})
$$ \bar{\mu}_{\min} (\shS)=\liminf_{q=p^{e}}  \mu_{\min} (\shS^q)/q \, .$$
This limit exists, since there exists for some number $k$ a
surjection $\oplus_j \O( \beta_j) \ra \shS(k) $ such that all $
\beta_j$ are positive. Then $\shS (k)$ is a quotient of an ample
bundle and so all its quotients have positive degree. This holds
also for all its Frobenius pull-backs, hence  $\mu_{\min}
((\shS(k))^q) \geq 0$ and the limit is $\geq 0$. Thus $\mu_{\min}
(\shS^q) \geq -q k \deg (\O_Y(1))$ for all $q$. Moreover, a theorem
of Langer implies that this limit is even a rational number, see
\cite{langersemistable}. The sheaf $\shS$ is called strongly
semistable if $\mu(\shS)=\bar{\mu}_{\min} (\shS)$; equivalently, if
all Frobenius pull-backs $\shS^q$ are semistable.

The degree of the variety $Y$ (with respect to $H$) is by definition
the top self intersection number $\deg(Y) = \deg(\O_Y(1))=H^t$. In
the following we will impose on a polarized variety $(Y, \O_Y(1))$
of dimension $t$ the condition that $H^i(Y,\O_Y(m))=0 $ for $i=1
\komdots t-1$ and all $m$. If $Y= \Proj R$, where $R$ is a
standard-graded Cohen Macaulay ring, this property holds true due to
\cite[Theorem 3.5.7]{brunsherzog}.

\begin{proposition}
\label{sheafproposition} Let $Y$ denote a normal projective variety
of dimension $t \geq 1$ over an algebraically closed field $K$ of
positive characteristic $p$. Let $\O_Y(1)$ denote a very ample
invertible sheaf on $Y$ such that $H^i(Y,\O(m))=0 $ for $i=1
\komdots t-1$. Suppose that the dualizing sheaf $\omega_Y$ on $Y$ is
invertible. Let $\shS$ denote a torsion-free coherent sheaf on $Y$.
Suppose that the stalk $\shS_y$ is free for every non-smooth point
$y \in Y$. Let
$$ \cdots \lra \shG_3 \lra \shG_2 \lra \shS \lra 0 $$
denote an exact complex of sheaves, where $\shG_j$ has type $\shG_j=
\bigoplus_{(k,j) } \O_Y(- \alpha_{k,j})$. Set $\shS_j=\im(\shG_{j+1}
\ra \shG_{j})= \ker( \shG_{j} \ra \shG_{j-1})$, $j \geq 2$, and
$\shS_1=\shS$. Fix $i=1 \komdots t$. Then for
$$m > - q \frac{  \bar{\mu}_{\min}( \shS_{t-i+1}) }{\deg(Y)} + \frac{\deg(\omega_Y)}{\deg(Y)} $$
we have $H^i (Y, \shS^q (m)) =0$.
\end{proposition}
\proof Note first that the Frobenius acts flat on the exact complex
and on the corresponding short exact sequences $0 \ra \shS_{j+1} \ra
\shG_{j+1} \ra \shS_j \ra 0$. This can be checked locally and is
true for the smooth points of $Y$. Over a singular point $y \in Y$
the sheaf $\shS$ is free, so these short exact sequences split
locally in a neighborhood of such a point and hence all the $\shS_j$
are also free in $y$. So also in these points the Frobenius
preserves the exactness of the complex.

Due to our assumption on $\O_Y(1)$ we have $H^i(Y,\shG_j(m))=0 $ for
$i=1 \komdots t-1$ and all $m$ and all $j \geq 2$. Hence from the
short exact sequences $0 \ra \shS_{j+1}(m) \ra \shG_{j+1}(m) \ra
\shS_j(m) \ra 0$ we can infer that
\begin{eqnarray*}
& & H^i (Y, \shS_j(m)) \cong H^{i+1}(Y, \shS_{j+1}(m)) \, \, \,
\mbox{ isomorphisms  for  } i=1 \komdots \dimty-2, \cr & & H^{\dimty
-1} (Y, \shS_j(m)) \subseteq H^{ \dimty}(Y, \shS_{j+1}(m)) \, \, \,
\mbox{ injection for } t \geq 2 , \cr & & H^{\dimty} (Y, \shG_{j+1}
(m)) \ra H^{\dimty} (Y, \shS_j(m)) \, \, \, \mbox{ surjection}.
\end{eqnarray*}
The same is true if we replace $S_j$ and $G_j$ by their Frobenius
pull-backs $S_j^q$ and $G_ j^q$. For $i=1 \komdots \dimty$ we find
$$H^{i}(Y, \shS^q_1 (m)) \!\cong \! H^{i+1}(Y, \shS^q_2 (m)) \! \congdots \!
H^{\dimty -1} (Y, \shS^q_{ \dimty-i} (m)) \! \subseteq \! H^{\dimty} (Y, \shS^q _{\dimty-i+1}(m)) .$$
So we only have to look at $H^\dimty (Y, \shS^q_{t-i+1}(m))$,
which is by Serre duality dual to
$\Hom(\shS^q_{t-i+1}(m), \omega_Y)$, see \cite[Theorem III.7.6]{hartshornealgebraic}.
Suppose now that $m$ fulfills the numerical condition.
Then
\begin{eqnarray*}
\mu_{\min}(\shS^q_{t-i+1}(m)) &=& \mu_{\min}(\shS^q_{t-i+1})+ m \deg(Y) \cr
&\geq & q \bar{\mu}_{\min} (\shS_{t-i+1}) + m \deg (Y) \cr
&>& q  \bar{\mu}_{\min}( \shS_{t-i+1})
+ \big(- q \frac{ \bar{\mu}_{\min}( \shS_{t-i+1}) }{\deg(Y)}
+ \frac{\deg(\omega_Y)}{\deg(Y)} \big) \deg(Y) \cr
&=& \deg(\omega_Y)  \, .
\end{eqnarray*}
So for these $m$ there are no non-trivial mappings from
$\shS_{t-i+1}^q(m)$ to $\omega_Y$ and therefore $H^{t}(Y, \shS_{t-i+1}^q(m))=0$.
\qed

\begin{remark}
The dualizing sheaf $\omega_Y$ on the projective variety $Y
\subseteq \PP^N$ is invertible under the condition that $Y$ is
locally a complete intersection in $\PP^N$ and in particular if $Y$
is smooth (see \cite[Theorem III.7.11 and Corollary
III.7.12]{hartshornealgebraic}. If $\omega_Y$ is not invertible, but
torsion-free, then we may replace $\deg(\omega_Y)$ by
$\mu_{\max}(\omega_Y)$ to get the same statement as in Proposition
\ref{sheafproposition}.
\end{remark}

\section{An inclusion bound for tight closure}
\label{inclusion}

We first fix the following situation, with which we will deal in
this section.

\begin{situation}
\label{situation} Let $K$ denote an algebraically closed field of
characteristic $p >0$. Let $R$ denote a standard-graded normal
Cohen-Macaulay domain of dimension $t+1 \geq 2$ over $K$ with
corresponding projective normal variety $Y=\Proj R$. Suppose that
the dualizing sheaf $\omega_Y$ of $Y$ is invertible. Let $I
\subseteq R$ denote a homogeneous $R_+$-primary ideal. Let
$$ \cdots  \lra F_2= \bigoplus_{(k,2)}  R(- \alpha_{k,2})
\lra F_1= \bigoplus_{(k,1)} R(- \alpha_{k,1}) \lra I \lra 0 \, ,$$
denote a homogeneous complex of graded $R$-modules which is exact on $D(R_+)$.
Let
$$ \cdots \lra \shG_2=\bigoplus_{(k,2)} \O(- \alpha_{(k,2)}) \lra
\shG_1= \bigoplus_{(k,1)} \O(- \alpha_{k,1})
 \lra \O_Y \lra 0 \, $$
denote the corresponding exact complex of sheaves on $Y$. Denote by
$\Syz_j = \ker (\shG_j \ra \shG_{j-1})$ the locally free kernel
sheaves on $Y$, and set $\Syz_j(m)= \Syz_j \otimes \O_Y(m)$. Let
$\nu =- \bar{\mu}_{\min} (\Syz_t) /\deg (Y)$, where $t$ is the
dimension of $Y$.
\end{situation}

\begin{theorem}
\label{theoreminclusion} Suppose the situation and notation
described in \ref{situation}. Then for all prime powers $q=p^{e}$ we
have the inclusion $R_{> q \nu + \frac{\deg(\omega_Y) }{\deg(Y)}}
\subseteq I^{[q]}$.
\end{theorem}
\proof Since $I$ is primary all the syzygy sheaves occurring in the
resolution on $Y$ are locally free and hence we may apply
Proposition \ref{sheafproposition}. Fix a prime power $q=p^{e}$. Let
$h \in R$ denote a homogeneous element of degree $m > q \nu +
\frac{\deg(\omega_Y) }{\deg(Y)}$. This gives via the short exact
sequence on $Y$,
$$0 \lra \Syz(f_1^q \komdots f_n^q)(m) \lra \bigoplus_{i=1}^n \O_Y(m -qd_i)
\stackrel{f_1^q \komdots f_n^q}{\lra} \O_Y(m) \lra 0 $$ rise to a
cohomology class $\delta(h) \in H^1(Y,\Syz(f_1^q \komdots
f_n^q)(m))$, where $$\Syz(f_1^q \komdots f_n^q)(m) =(F^{e*}(\Syz(f_1
\komdots f_n))) (m)= \shS^q (m) \, ,$$ $\shS = \shS_1= \Syz(f_1
\komdots f_n)$. It is enough to show that $\delta (h)=0$, for then
$h \in I^{[q]} \Gamma(D(R_+), \O) =I^{[q]}$, since $R$ is normal.
But this follows from Proposition \ref{sheafproposition} applied to
$\shS= \Syz(f_1 \komdots f_n)$ and $i=1$. \qed

\begin{remark}
\label{resolutionremark} We do not insist that the ``resolution'' of
the ideal is exact on the whole $\Spec R$ nor that it is minimal,
but it is likely that a minimal resolution will give us in general a
better bound $\nu$. For example we can always use the Koszul complex
given by ideal generators of the $R_+$-primary ideal $I$.
\end{remark}

The next theorem gives an inclusion bound for tight closure. Recall
that the tight closure of an ideal $I \subseteq R$ in a noetherian
domain containing a field of positive characteristic $p$ is by
definition the ideal
$$I^* = \{f \in R: \exists 0 \neq  c \in R \mbox{ such that } cf^q \in I^{[q]}
\mbox{ for all } q=p^{e} \} \, .$$
See \cite{hunekeapplication} for basic properties of this closure operation.

\begin{theorem}
\label{tightinclusion} Suppose the situation described in
\ref{situation}. Then we have the inclusion $R_{\geq  \nu} \subseteq
I^*$.
\end{theorem}
\proof Let $f \in R$ be a homogeneous element of degree $\deg (f)=m
\geq \nu =- \bar{\mu}_{\min}(\Syz_t) /\deg(Y)$. Due to the
definition of tight closure we have to show that $cf^q \in I^{[q]}$
holds for some $c \neq 0$ and all prime powers $q$. Let $c \neq 0$
be any homogeneous element of degree $ >  \deg (\omega_Y)/\deg (Y)$.
Then $\deg (cf^q) = qm +\deg(c) > q \nu + \deg(\omega_Y)/\deg (Y)$
and therefore $cf^q \in I^{[q]}$ by Theorem \ref{theoreminclusion}.
\qed

\begin{remark}
Suppose that $R$ fulfills the condition of the situation described
in \ref{situation} and let $I=(f_1 \komdots f_n)$ denote an ideal
generated by a full regular system of homogeneous parameters of
degree $\deg(f_i)=d_i$ (so $n=t+1$). Then the Koszul resolution of
these elements gives a resolution on $Y=\Proj R$ such that the
top-dimensional syzygy bundle is invertible, namely
$$\Syz_t (m)=\shG_{t+1}(m)= \O_Y(m-d_1 - \cdots -d_{t+1} ) \, $$
Then Theorem \ref{tightinclusion} gives the known (even without the
condition Cohen-Macaulay) inclusion bound $R_{\geq d_1 \plusdots
d_n} \subseteq (f_1 \komdots f_n)^*$, see \cite[Theorem
2.9]{hunekeparameter}.

The next easiest case is then the $R_+$-primary homogeneous ideal
$I$ has finite projective dimension (it is again enough to impose
the exactness only on $D(R_+)$). In this case the resolution on $Y$
looks like
$$0 \lra \shG_{t+1} \lra \shG_t \lradots  \shG_1 \lra \O_Y\lra 0$$
and the top-dimensional syzygy bundle is $\Syz_t = \shG_{t+1}
=\bigoplus_k \O_Y( - \alpha_{k,t+1})$, and therefore
$$\mu_{\min}(\Syz_t)= \deg(Y) \min_k \{ - \alpha_{k,t+1}\} = -
\deg(Y) \max_k \{ \alpha_{k,t+1} \} \, .$$ The corresponding
inclusion bound was proved in \cite[Theorem
5.11]{hunekesmithkodaira}. Such a situation arises for example if
$I$ is generated by a set of monomials in a system of homogeneous
parameters.
\end{remark}

The following easy corollary unifies two known inclusion bounds for
tight closure given by K. Smith (see \cite[Propositions 3.1 and
3.3]{smithgraded}), namely that $R_{\geq \sum_{i=1}^n \deg (f_i)}
\subseteq I^*$ and that $R_{\geq \dim(R) \max _i \{ \deg(f_i)\}}
\subseteq I^*$.

\begin{corollary}
\label{tightcorollary} Suppose the situation described in
\ref{situation} and suppose that the homogeneous $R_+$-primary ideal
$I=(f_1 \komdots f_n)$ is generated by homogeneous elements of
degree $d_i= \deg(f_i)$. Set $d= \max _{1 \leq i_1 < \ldots <
i_{\dim (R)} \leq n}( d_{i_1} \plusdots d_{i_{\dim(R)}} )$. Then $R
_{\geq d} \subseteq I^*$.
\end{corollary}
\proof
We consider the Koszul resolution of $I=(f_1 \komdots f_n)$,
which is exact outside the origin.
This gives the surjection
$$\bigoplus _{1 \leq i_1 < \ldots < i_{\dim (R)} \leq n}
\O(-d_{i_1} \minusdots  d_{i_{\dim(R)}} ) \lra  \Syz_{\dim(R)-1 } \lra 0 $$
which shows that
\begin{eqnarray*}
\bar{\mu}_{\min} (\Syz_{\dim(R)-1}) &\geq&  \bar{\mu}_{\min}
\big(\bigoplus _{1 \leq i_1 < \ldots < i_{\dim (R)} \leq n}
\O(-d_{i_1} \minusdots d_{i_{\dim(R)}}) \big) \cr &=& - \max _{1
\leq i_1 < \ldots < i_{\dim (R)} \leq n} \{ d_{i_1} \plusdots
d_{i_{\dim(R)}}   \} \deg (Y) \, .
\end{eqnarray*}
Hence $\nu =- \bar{\mu}_{\min} ( \Syz_{\dim(R)-1})/ \deg (Y) \leq
\max \{d_{i_1} \plusdots d_{i_{\dim(R)}} \}$ and Theorem
\ref{tightinclusion} applies. \qed

\begin{remark}
If the dimension of $R$ is two, then Theorem \ref{tightinclusion}
was proved in \cite[Theorem 6.4]{brennerslope} using somewhat more
geometric methods. In this case $Y=  \Proj R$ is a smooth projective
curve and the top syzygy bundle is just the first syzygy bundle, and
the result also holds in characteristic zero for solid closure. See
\cite{brennerslope} and \cite{brennercomputationtight} for concrete
computations of the number $\nu$ in this case. It is in general
difficult to compute the number $\nu$ of the theorem, as it is
difficult to compute the minimal slope of a locally free sheaf.
\end{remark}

The following corollary gives an inclusion bound for tight closure
under the condition that the top-dimensional syzygy bundle is
strongly semistable. In the two-dimensional situation this bound is
exact, in the sense that below this bound an element belongs to the
tight closure only if it belongs to the ideal itself, see
\cite[Theorem 8.4]{brennerslope}.

\begin{corollary}
\label{topstable} Suppose the situation described in \ref{situation}
and let $I=(f_1 \komdots f_n)$ be generated by homogeneous elements
of degree $d_i=\deg(f_i)$. Let $F_\bullet \ra I$ denote the Koszul
complex and suppose that the top-dimensional syzygy bundle $\Syz_t$
is strongly semistable. Set $d= (\dim (R)-1) (d_1 \plusdots d_n)/
(n-1)$. Then $R_{\geq d} \subseteq I^*$.
\end{corollary}
\proof
The condition strongly semistable means that $\mu(\Syz_t)= \bar{\mu}_{\min}(\Syz_t)$.
So we only have to compute the degree and the rank of $\Syz_t$.
It is easy to compute that $\det(\Syz_t)= \O_Y(\binom{n-2}{t-1}(- \sum_{i=1}^n d_i))$,
hence
$$\deg( \Syz_t) = \binom{n-2}{t-1}(- \sum_{i=1}^n d_i) \deg (Y)$$
and $\rk (\Syz_t)= \binom{n-1}{t}$.
Therefore
$$\mu( \Syz_t)=  \binom{n-2}{t-1} (- \sum_{i=1}^n d_i) \deg (Y) / \binom{n-1}{t}
= \frac{t}{n-1} (- \sum_{i=1}^n d_i) \deg (Y) $$
and $\nu = \frac{t}{n-1} ( \sum_{i=1}^n d_i)$.
\qed

\begin{remark}
As the proofs of Theorem \ref{tightinclusion} and Proposition
\ref{sheafproposition} show, Corollary \ref{topstable} is also true
under the weaker condition that there does not exist any non-trivial
mapping $\Syz_t^q \ra \shL$ to any invertible sheaf $\shL$
contradicting the semistability of $\Syz_t^q$ for all $q=p^{e}$.
\end{remark}

\begin{example}
Theorem \ref{tightinclusion} applies in particular when $R$ is a
normal complete intersection domain. Let $R=K[X_1 \komdots X_N]/(H_1
\komdots H_r)$, where $H_j$ are homogeneous forms of degree
$\delta_j$. Then $\omega_Y = \O(\sum_j \delta_j -N)$. Therefore the
number $ \deg(\omega_Y)/ \deg(Y)= \sum_j \delta_j -N$ is just the
$a$-invariant of $R$.
\end{example}

\begin{example}
We want to apply Corollary \ref{topstable} to the computation of the
tight closure $(x^a,y^a,z^a,w^a)^*$ in $R=K[x,y,z,w]/(H)$, where $H$
is supposed to be a polynomial of degree $4$ defining a smooth
projective (hyper-)surface
$$Y= V_+(H)=\Proj R \subset \PP^3 =\Proj K[x,y,z,w]$$
of degree $4$; hence $Y$ is a $K 3$ surface. Our result will only
hold true for generic choice of $H$. We look at the Koszul complex
on $\PP^3$ defined by $x^a,y^a,z^a,w^a$ and break it up to get
$$0 \lra \Syz_2 \cong \bigwedge^2 \Syz \lra \bigoplus_6 \O_{\PP^3}(-2a)
\lra \bigoplus_4 \O_{\PP^3}(-a) \lra \O_{\PP^3} \lra 0 \, .$$
Suppose first that $K$ is an algebraically closed field of
characteristic $0$. It is easy to see that the syzygy bundle
$\Syz=\Syz(x^a,y^a,z^a,w^a)$ is semistable on $\PP^3$
\cite[Corollary 3.6 or Corollary 6.4]{brennerlookingstable}.
Therefore also the exterior power $\Syz_2 \cong \bigwedge^2 \Syz$ is
semistable on $\PP^3$. By the restriction theorem of Flenner
\cite[Theorem 1.2]{flennerrestriction} it follows that the
restriction $\Syz_2 \!|_Y$ is also semistable on the generic
hypersurface $Y=V_+(H)$.

On the other hand, due to the Theorem of Noether (see \cite[\S
IV.4]{haramp}), every curve on the generic surface of degree $4$ in
$\PP^3$ is a complete intersection and $R=K[x,y,z,w]/(H)$ is a
factorial domain for generic $H$ of degree $4$. It follows that the
cotangent bundle $\Omega_Y$ on $Y=V_+(H)$ is semistable. For the
semistability of a rank two bundle we only have to look at mappings
$\shL  \ra \Omega_Y$, where $\shL$ is invertible. But since $\shL=
\O_Y(k)$, the semistability follows, since $Y$ is a $K3$ surface and
so $\Omega_Y$ has degree $0$ but does not have any global
non-trivial section (see \cite[IV. 5]{griffithsharris}).

So for $H$ generic the relevant second syzygy bundle $\Syz_2 \!|_Y$
and the cotangent bundle $\Omega_Y$ are both semistable in
characteristic $0$. Since the $\QQ$-rational points are dense in
$\AA^N_K$, there exist also such polynomials $H$ with rational
coefficients and then also with integer coefficients. We consider
such a polynomial $H$ with integer coefficients as defining a family
of quartics over $\Spec \ZZ$. Since semistability is an open
property, we infer that the second syzygy bundle and the cotangent
bundle are also semistable on $Y_p=V_+(H_p)$ for $p \gg 0$.

By the semistability of $\Omega_{Y_p}$ ($p \gg 0$), the maximal
slope of $\Omega_{Y_p}$ is $\leq 0$. A theorem of Langer
\cite[Corollary 2.4 and Corollary 6.3]{langersemistable} shows then
that every semistable bundle on $Y_p$ is already strongly
semistable. Hence the second syzygy bundle is also strongly
semistable. Therefore we are in the situation of Corollary
\ref{topstable} and we compute $d=8a/3$. Thus
$$ R_{8a/3}\subseteq (x^a,y^a,z^a,w^a)^*$$ holds in $R=K[x,y,z,w]/(H)$ for $H$
generic of degree $4$ and for $p \gg 0$. The first non-trivial
instance is for $a=3$. In fact for the (non-generic) Fermat quartic
$x^4+y^4+z^4+w^4=0$ it was proved by Singh in \cite[Theorem
4.1]{singhcomputation} directly that $x^2y^2z^2w^2 \in
(x^3,y^3,z^3,w^3)^*$.
\end{example}

For the next corollary we recall the definition of the Frobenius
closure. Suppose that $R$ is a noetherian ring containing a field of
positive characteristic $p >0$, and let $I$ denote an ideal. Then
the Frobenius closure of $I$ is defined by
$$I^F =\{ f \in R:\, \exists q=p^{e} \mbox{ such that } f^q \in I^{[q]} \} \, .$$
It is easy to see that the Frobenius closure of an ideal is contained in its
tight closure.

\begin{corollary}
\label{frobenius} Suppose the situation described in
\ref{situation}. Then $R_{> \nu} \subseteq I^F$, the Frobenius
closure of $I$.
\end{corollary}
\proof Let $f$ denote a homogeneous element of degree $m=\deg(f) >
\nu=- \bar{\mu}_{\min} (\Syz_t)/\deg(Y)$. Then we just have to take
a prime power $q=p^{e}$ such that $\deg(f^q)=qm >q \nu+
\deg(\omega_Y)/ \deg(Y)$ holds. Then $f^q \in I^{[q]}$ holds due to
Theorem \ref{theoreminclusion}. \qed

\begin{example}
Corollary \ref{frobenius} is not true for $R_{\geq \nu}$ instead of
$R_{> \nu }$. This is already clear for parameter ideals in
dimension two, say for $(x,y)$ in $R=K[x,y,z]/(H)$, where $H$
defines a smooth projective curve $Y= \Proj R =V_+(H) \subset
\PP^2$. Here we have the resolution
$$0 \lra \O_Y(-2) \cong \Syz(x,y)(0) \lra \O_Y(-1) \oplus \O_Y(-1) \stackrel{x,y}{\lra}
\O_Y \lra 0 \, .$$ Hence we get $\nu =2$, but an element of degree
two (say $z^2$) does not in general belong to the Frobenius closure
of $(x,y)$.
\end{example}

\begin{remark}
A problem of Katzman and Sharp (see \cite{katzmansharpfrobenius})
asks in its strongest form: does there exist a number $b$ such that
whenever $f \in I^F$ holds, then already $f^{p^b} \in I^{[p^b]}$
holds. A positive answer (together with the knowledge of a bound for
the number $b$) to this question would give a finite test to check
whether a given element $f$ belongs to the Frobenius closure $I^F$
or not. For those elements which belong to $I^{F}$ because of
Corollary \ref{frobenius} (due to degree reasons, so to say), the
answer is yes, at least in the sense that for $f$ fulfilling
$\deg(f) \geq \nu + \epsilon$ ($\epsilon
>0$) we have $ \deg (f^q) = q \deg(f) \geq q \nu + q \epsilon$, so
the condition $q \epsilon > \deg(\omega_Y)/\deg(Y)$ is sufficient to
ensure that $f^q \in I^{[q]}$. It is however possible that elements
of degree $\deg(f) \leq \nu$ belong to the Frobenius closure.
\end{remark}

\section{The Castelnuovo-Mumford regularity of Frobenius powers}
\label{cmregularity}

We recall briefly the notion of the Castelnuovo-Mumford regularity following
\cite[Definition 15.2.9]{brodmannsharp}.
Let $R$ denote a standard-graded ring and let $M$ denote a finitely generated
graded $R$-module.
Then the Castelnuovo-Mumford regularity of $M$ (or regularity of $M$ for short)
is
$$ \reg(M) = \sup \{ \engrad (H^i_{R_+} (M)) +i :\, 0 \leq i \leq \dim M \} \, ,$$
where $\engrad (N)$ of a graded $R$-module $N$ denotes the maximal
degree $e$ such that $N_e \neq 0$. For a number $l$ we define the
regularity $\reg^l (M)$ at and above level $l$ by
$$\reg^l(M)= \sup \{ \engrad (H^i_{R_+} (M))+i :\, l \leq i \leq \dim M \} \, ,$$

A question of M. Katzman raised in
\cite[Introduction]{katzmanfrobenius} asks how the regularity of the
Frobenius powers $I^{[q]}$ behaves, in particular whether there
exists a linear bound $\reg(I^{[q]}) \leq C_1 q + C_0$. Such a
linear bound for the regularity of the Frobenius powers of an ideal
was recently given by M. Chardin in \cite[Theorem
2.3]{chardinregularitypowers}. The following theorem gives a better
linear bound for the regularity of Frobenius powers of $I$ in terms
of the slope of the syzygy bundles.

\begin{theorem}
\label{regularitybound} Let $K$ denote an algebraically closed field
of positive characteristic $p$. Let $R$ denote a standard-graded
normal Cohen-Macaulay $K$-domain of dimension $t+1 \geq 2$. Let
$I=(f_1 \komdots f_n) \subseteq R$ denote a homogeneous ideal
generated by homogeneous elements of degree $d_i =\deg(f_i)$.
Suppose that the dualizing sheaf $\omega_Y$ on $Y=\Proj R$ is
invertible. Suppose that the points $y \in \sup (\O_Y/\shI)$ are
smooth points of $Y$. Let $F_\bullet \ra I$ denote a graded free
resolution with corresponding exact complex of sheaves on $Y$, $
\shG_\bullet \ra \shI \subseteq \O_Y $. Set $\Syz_j = \ker(\shG_j
\ra \shG_{j-1})$. Then we have for the Castelnuovo-Mumford
regularity of the Frobenius powers $I^{[q]}$ the linear bound $\reg
(I^{[q]}) \leq C_1q + C_0 $, where
$$C_1= \max \{ d_i, i=1 \komdots n, \,
- \frac{ \bar{\mu}_{\min}(\Syz_j)}{\deg(Y)},\, j=1 \komdots t =\dim(Y) \} \mbox{ and }   $$
$$ C_0 = \max \{ \reg(R), \frac{\deg(\omega_Y)}{\deg(Y) } \}
\, .$$
\end{theorem}
\proof
The ideal generators define for $q=p^{e}$ the homogeneous short exact sequences
$$ 0 \lra \Syz(f_1^q \komdots f_n^q) \lra \bigoplus_{i=1}^n R(-qd_i)
\stackrel{f_1^q \komdots f_n^q}{ \lra} I^{[q]} \lra 0$$ of graded
$R$-modules. It is an easy exercise \cite[Exc.
15.2.15]{brodmannsharp} to show that for a short exact sequence $0
\ra L \ra M \ra N \ra 0$ we have $\reg(N) \leq \max \{ \reg^1(L) -1,
\reg (M) \}$. We have $\reg (R(-qd)) = \reg (R) +qd$ and
$$\reg ( \bigoplus_{i=1}^n R(-qd_i)) = \max_i \{\reg (R(-qd_i)) \}
=\reg (R) + q \max_i \{ d_i \} \, ,$$ which gives the first terms in
the definition of $C_1$ and $C_0$ respectively. Hence it is enough
to give a linear bound for $\reg^1(\Syz(f_1^q \komdots f_n^q))$.
Moreover, the long exact local cohomology sequence associated to the
above short exact sequence gives
$$\lra H^0_{R_+}(I^{[q]}) \lra
H^1_{R_+} (\Syz(f_1^q \komdots f_n^q)) \lra \bigoplus_{i=1}^n
H^1_{R_+} (R(-qd_i)) \lra \, . $$ The term on the right is $0$,
since $R$ is Cohen-Macaulay, and the term on the left is $0$, since
$R$ is a domain. Therefore $H^1_{R_+} (\Syz(f_1^q \komdots f_n^q))
=0$ and we have to find a linear bound for $ \reg^2(\Syz(f_1^q
\komdots f_n^q)) =\reg^1(\Syz(f_1^q \komdots f_n^q))$. We have
$H^i_{R_+}(\Syz(f_1^q \komdots f_n^q) )= H^{i-1} (D(R_+), \Syz(f_1^q
\komdots f_n^q) \widetilde{\, } \, )$ for $i \geq 2$ due to the long
exact sequence relating local cohomology with sheaf cohomology.
Denote now by $\Syz (f_1^q \komdots f_n^q) $ the corresponding
torsion-free sheaf on $Y=\Proj R$. On $Y$ we have the short exact
sequences of sheaves
$$0 \lra \Syz (f_1^q \komdots f_n^q)
\lra \bigoplus_{i=1}^n \O_Y(-qd_i) \lra \shI^{[q]} \lra 0  \, .$$
We may compute the cohomology  as
$$H^{i}(D_+(R), \Syz(f_1^q \komdots f_n^q) \widetilde{\,} \,  )_m
= H^i(Y, \Syz(f_1^q \komdots f_n^q) (m)) \, .$$ Note that the syzygy
bundle $\Syz(f_1 \komdots f_n)$ is free by assumption in the
singular points of $Y$. Hence we are in the situation of Proposition
\ref{sheafproposition} with $\shS=\Syz(f_1 \komdots f_n)$; therefore
$H^{i} (Y, \Syz(f_1^q \komdots f_n^q)(m))=0$ ($i=1 \komdots t$)
holds for $m> \max_{j=1 \komdots t} \{- q \frac{ \bar{\mu}_{\min}
(\Syz_j)}{ \deg(Y)} \} + \frac{\deg(\omega_Y)}{\deg(Y)}$, which
proves the theorem. \qed

\begin{remark}
The Castelnuovo-Mumford regularity of a standard-graded Cohen-Macaulay domain $R$
is just $\reg(R)= \engrad (H^{\dim ( R)}_{R_+} (R)) + \dim (R)$.
The end of the top-dimensional local cohomology module of a graded ring is also called its
$a$-invariant, see \cite[13.4.7]{brodmannsharp}, hence $\reg(R)=a + \dim (R)$.
If $R$ is Gorenstein, then $R(a)$ is the canonical module of $R$ and
$\omega_Y=\O_Y(a)$ is the dualizing sheaf on $Y=\Proj R$. So in this case
the quotient $\deg (\omega_Y) / \deg (Y)= a \deg (Y)/ \deg (Y) =a$
equals also the $a$-invariant.
\end{remark}

\begin{remark}
The surjection $\bigoplus_{(k,j+1)} \O_Y(- \alpha_{k,j+1}) \ra
\Syz_j \ra 0$ gives at once the bound $\bar{\mu}_{\min} (\Syz_j)
\geq \bar{\mu}_{\min}(\bigoplus_{(k,j+1)} \O_Y(- \alpha_{k,j+1}))
\!= \! - \max \{\alpha_{k,j+1} \} \deg(Y)$. Therefore we get for the
constant $C_1$ coming from Theorem \ref{regularitybound} the
estimate $C_1 \leq \max \{ \alpha_{k,j}:\, j=1 \komdots t+1 = \dim
(R) \,  \} =C_1'$. This number $C_1'$ is the coefficient for the
linear bound which M. Chardin has obtained in \cite[Theorem
2.3]{chardinregularitypowers}. This bound corresponds to the
inclusion bounds for tight closure of K. Smith which we obtained in
Corollary \ref{tightcorollary}. The following standard example of
tight closure theory shows already the difference between the
Chardin-Smith bound and the slope bound.
\end{remark}

\begin{example}
Consider the ideal $I=(x^2,y^2,z^2)$ in $R=K[z,y,z]/(x^3+y^3+z^3)$,
${\rm char} (K) \neq 3$. We compute the bound coming from Theorem
\ref{regularitybound} for the regularity of the Frobenius powers
$I^{[q]}=(x^{2q}, y^{2q},z^{2q})$. We first observe that we may
consider the curve equation $0=x^3+y^3+z^3=xx^2+yy^2+zz^2$ as a
global section of the syzygy bundle of degree $3$. Since this
section has no zero on $Y= \Proj R$, we get the short exact sequence
$$0 \lra \O_Y \lra \Syz(x^2,y^2,z^2)(3) \lra \O_Y \lra0 \, .$$
This shows that the syzygy bundle is strongly semistable and
therefore $\bar{\mu}_{\min} ( \Syz(x^2,y^2,z^2)(0))= -6 \deg(Y)/2=
-9 $. So $C_1=3$ and we get altogether the bound $\reg( I^{[q]})
\leq 3q+2$.

Since $\, \Syz(x^2,y^2,z^2)(3)\, $ is not generated by its global
sections, because the section just mentioned is the only section.
Hence a surjection $ \bigoplus_k \O(-\alpha_k) \ra
\Syz(x^2,y^2,z^2)(0)$ is only possible for $ \max _k \{ \alpha_k \}
\geq 4$. So the linear bound for the regularity which you get by
considering only the degrees in a resolution is worse than the slope
bound.
\end{example}

\bibliographystyle{plain}

\bibliography{bibliothek}

\end{document}